\documentclass[12pt]{article}
\usepackage{amsfonts}
\usepackage{amsmath}

\newcommand{\BEAS}{\begin{eqnarray*}}
\newcommand{\EEAS}{\end{eqnarray*}}
\newcommand{\BEA}{\begin{eqnarray}}
\newcommand{\EEA}{\end{eqnarray}}
\newcommand{\BEQ}{\begin{equation}}
\newcommand{\EEQ}{\end{equation}}
\newcommand{\BIT}{\begin{itemize}}
\newcommand{\EIT}{\end{itemize}}
\newcommand{\BNUM}{\begin{enumerate}}
\newcommand{\ENUM}{\end{enumerate}}
\newcommand{\BA}{\begin{array}}
\newcommand{\EA}{\end{array}}
\newtheorem{theorem}{Theorem}

\newtheorem{proposition}[theorem]{Proposition}

\newenvironment{proof}{\begin{quote}\begin{small}\textbf{Proof.}}{\end{small}\QED\end{quote}}
\newcommand{\reals}{{\mbox{\bf R}}}

\newcommand{\symm}{{\mbox{\bf S}}}  

\newcommand{\Rank}{\mathop{\bf Rank}}

\newcommand{\Card}{\mathop{\bf Card}}

\newcommand{\Co}{{\mathop {\bf Co}}}

\newcommand{\QED}{~~\rule[-1pt]{6pt}{6pt}}

\begin{document}

\title{A semidefinite representation for some minimum cardinality
problems}
\author{Alexandre d'Aspremont\thanks{%
Information Systems Laboratory \& Department of Management Science
and Engineering, Terman bldg., Stanford University, Stanford, CA
94305-4023, USA. Internet: aspremon@stanford.edu}} \maketitle

\begin{abstract}
Using techniques developed in \cite{Lass02a}, we show that some minimum
cardinality problems subject to linear inequalities can be represented as
finite sequences of semidefinite programs. In particular, we provide a
semidefinite representation of the minimum rank problem on positive
semidefinite matrices. We also use this technique to cast the problem of
finding convex lower bounds on the objective as a semidefinite program.

\end{abstract}

\subsection*{Notation}

We note $\reals[x_{1},...,x_{n}]$ (or $\reals[x]$ when there is no
ambiguity) the ring of multivariate polynomials
$p(x)=p(x_{1},...,x_{n})$ on
a variable $x\in \reals^{n}$. We say that $p(x)\in \reals[x]$ is $SOS$ when $%
p(x)$ is a sum of squares of polynomials in $\reals[x]$. For $x\in \reals%
^{n} $, $\Card(x)$ will be the cardinal of the set $\{i:x_{i}\neq
0\}$. We note $\symm^{n}$ the set of $n\times n$ symmetric
matrices. For multivariate polynomials, we adopt the multiindex
notation $p(x)=\sum_{\alpha }p_{\alpha }x^{\alpha }$, where
$x^{\alpha }:=x_{1}^{\alpha _{1}}x_{2}^{\alpha
_{2}}...x_{n}^{\alpha _{n}}$, and we note $d=\sum_{i=1}^{n}\alpha
_{i}$ the degree of $p(x).$ $\reals_{d}[x]$ is the set of
polynomials of degree at most $d$. Finally, $C(p)$ will be the
\textit{Newton polytope} of the polynomial $p(x)$, with
$C(p)=\Co\left( \left\{ \alpha :p_{\alpha }\neq 0\right\} \right)
$.

\section{Introduction}

Given a convex set $\mathcal{C}\subset \reals^{n}$, we are
interested in solving the following problem:
\begin{equation}
\begin{array}{ll}
\text{minimize} & \Card(x) \\
\text{subject to} & x\in \mathcal{C},
\end{array}
\label{pb-min-card}
\end{equation}
in the particular case where $\mathcal{C}$ is described by a set
of linear inequalities. Except in certain rare instances, this
problem is very hard to solve (see \cite{Vand96}). Excellent
heuristics exist however, a classical
one (see \cite{Hass99} for example) replacing the function $\Card(x)$ by $%
\left\| x\right\| _{1}$, its largest convex lower bound on the
unit cube.

A related problem is that of minimizing the rank of a p.s.d.
matrix subject to LMI constraints:
\begin{equation}
\begin{array}{ll}
\text{minimize} & \Rank(X) \\
\text{subject to} & X\in \mathcal{C},
\end{array}
\label{pb-min-rank}
\end{equation}
where $\mathcal{C}$ is here an affine subset of the semidefinite
cone (a LMI). In this case also, minimizing the nuclear norm
$\left\| X\right\| _{\ast }$ of $X$ will produce excellent
approximate solutions (see \cite {Faze00}).

In this paper, using results by \cite{Cass84}, \cite{Shor87},
\cite{Puti93},
\cite{Choi95}, \cite{Nest00}, \cite{Lass01}, \cite{Pari01} and \cite{Lass02a}%
, we show that the $\mathbf{Min}\Card(x)$ and
$\mathbf{Min}\Rank(X)$ problems in (\ref{pb-min-card}) and
(\ref{pb-min-rank}) are equivalent to large semidefinite programs
(see \cite{Nest94}). To be precise, based on a
reformulation \`{a} la \cite{Shor87} of problems (\ref{pb-min-card}) and (%
\ref{pb-min-rank}), we use the technique in \cite{Lass02a} to
produce a finite (possibly exponential) sequence of increasingly
tighter semidefinite relaxations.

The rest of the paper is organized as follows. In section
\ref{s:sdr-sos}, we recall some key definitions and properties on
semidefinite representability and the sum of squares
representation of positive polynomials. We also summarize the
application of these representations to semialgebraic problems. In
section \ref{s:sdp-sap}, we show that both the
$\mathbf{Min}\Card(x)$ and the $\mathbf{Min}\Rank(X)$ problems are
equivalent to large scale semidefinite programs. Based on the work
by \cite {Puti93}, \cite{Nest00} and \cite{Lass02a} we explicitly
construct in section \ref{s:sdp-relax} a sequence of semidefinite
programs solving problems (\ref{pb-min-card}) and
(\ref{pb-min-rank}). We also show how the problem of finding
optimal convex lower bounds on the objective function can be
represented in a similar way. Finally, in section
\ref{s:complexity}, we discuss the complexity of these techniques.

\section{Sums of squares and semidefinite programming}

\label{s:sdr-sos}We quickly recall here some key definitions and
properties linking semidefinite and semialgebraic problems.

\label{ss:sdr-sos-sos}Hilbert's 17$^{th}$ problem (see
\cite{Rezn96} for an overview), which asked if all positive
polynomials could be written as sums of squares of other
polynomials, has a positive answer in dimension one. \cite{Nest00}
provides an efficient way of computing the $SOS$ representation of
a given positive univariate polynomial in the following result
from \cite{Nest00}.

\begin{proposition}
Let $p(x)\in \reals[x]$ be a univariate polynomial of degree $d$.
Then $p(x)$ for all $x\in \reals$ iff there exists a matrix $X\in
\symm^{v}$ such that:
\begin{equation}
p(x)=y_{x}^{T}Xy_{x}\text{, with }X\succeq 0,\quad \text{for all
}x\in \reals \label{prop-repres-univ-sos}
\end{equation}
with $v=\left\lceil d/2\right\rceil $ and $y_{x}=\left(
1,x,x^{2},...,x^{v}\right) $ is the list of univariate monomials
up to degree $v$.
\end{proposition}

The coefficients of the polynomials in the representation are then
computed as the eigenvectors of the matrix $X$.

In the general multivariate case, that representation property of
positive polynomials is lost. It can be shown (see \cite{Berg80})
that the set of multivariate $SOS$ polynomials is dense in the set
of positive polynomials, but there are simple examples of positive
polynomials that are not $SOS$. However, recent results in
semialgebraic geometry (see \cite{Cass84}, \cite{Schm91},
\cite{Puti93} or \cite{Puti99}) bridge the gap between positive
and $SOS$ polynomials on compact semialgebraic sets. We cite here
the result in \cite{Puti93}. Let $g_{k}(x)\in $
$\reals[x_{1},...,x_{n}]$ for $j=1,...,r,$ and we note $K$, the
semialgebraic set defined by $K=\left\{ x\in
\reals^{n}:g_{k}(x)\geq 0,k=1,...,r\right\} $. We suppose that $K$
is compact and that there exists $u(x)\in \reals[x_{1},...,x_{n}]$
such that $\left\{ u(x)\geq 0\right\} $ is compact with
\begin{equation}
u(x)=u_{0}(x)+\sum_{k=1}^{r}g_{k}(x)u_{k}(x),\quad \text{for all }x\in \reals%
^{n}  \label{assumpt-Putinar}
\end{equation}
where the polynomials $u_{k}(x)\in $ $\reals[x_{1},...,x_{n}]$ \
are $SOS$ for $k=1,...,r$. Under this assumption, we can represent
all polynomials positive on $K$ using $SOS$ polynomials as in
\cite{Puti93} or \cite{Puti99}.

\begin{proposition}
Suppose (\ref{assumpt-Putinar}) holds. A polynomial $p(x)\in $ $%
\reals%
[x_{1},...,x_{n}]$ is positive on $K$ iff:
\begin{equation}
p(x)=q_{0}(x)+\sum_{k=1}^{r}g_{k}(x)q_{k}(x),\quad \text{for all }x\in \reals%
^{n}  \label{prop-Putinar-sos}
\end{equation}
where the polynomials $q_{k}(x)\in \reals[x_{1},...,x_{n}]$\ are $SOS$ for $%
k=0,...,r$.
\end{proposition}

Now, as in \cite{Pari00} or \cite{Lass01}, we can write the
multivariate version of the relation (\ref {prop-repres-univ-sos})
mapping $SOS$ polynomials to the semidefinite cone.

\begin{proposition}
Let $p(x)\in $ $\reals[x_{1},...,x_{n}]$ be a polynomial and $K\
$a semialgebraic set defined by $K=\left\{ x\in
\reals^{n}:g_{k}(x)\geq
0,k=1,...,r\right\} $, satisfying assumption (\ref{assumpt-Putinar}). Then $%
p(x)\geq 0$ on $K$ iff there is an integer $m\in \mathbf{Z}_{+}$
and matrices $X_{k}\in \symm^{N}$, with $X_{k}\succeq 0$ for
$k=0,...,r$ such that:
\begin{equation}
p(x)=y_{x}^{T}X_{0}y_{x}+\sum_{k=1}^{r}\left(
y_{x}^{T}X_{k}y_{x}\right) g_{k}(x),\quad \text{for all }x\in
\reals^{n},  \label{prop-Lasserre-sos}
\end{equation}
where $N=\left\lceil \binom{n+m-1}{m}/2\right\rceil $ and
$y_{x}=\left(
1,x_{1},...,x_{n},x_{1}^{2},x_{1}x_{2},...,x_{1}x_{n},x_{2}x_{3},...,x_{n}^{2},...,x_{1}^{m},...,x_{n}^{m}\right)
$ is the vector of all monomials in $\reals[x_{1},...,x_{n}]$, up to degree $%
m$, listed in graded lexicographic order.
\end{proposition}

The result on polynomials above shows that testing the positivity
of a multivariate polynomial on a semialgebraic set $K$ satisfying
the assumption (\ref{assumpt-Putinar}) can be cast as a
semidefinite program. In general, the result in \cite{Lass01}
shows that all compact semialgebraic problems, i.e. problems
seeking to minimize a polynomial over a compact semialgebraic set,
are equivalent to large-scale semidefinite programs. This provides
a positive answer in the compact multivariate case to all the open
questions in \S 4.10.2 in \cite{Bent01}. A converse result is also
true (and much simpler). Because the positive semidefiniteness of
a matrix is equivalent to that of all its principal minors, all
semidefinite programs are semialgebraic programs, with additional
convexity and invariance properties.

The central result of moment theory exploited in \cite{Lass01}
sets polynomial positivity problems and moment problems as duals
(see e.g. \cite {Berg80}). Let $s$ be a positive semidefinite
sequence $s\in \reals^{N}$, we have
\[\BA{c}
s\text{ is p.s.d.}\\
\Updownarrow\\
\left\langle s,p_{\alpha}\right\rangle \geq 0,\quad \text{for all
}p(x)\in \reals_{m}[x]\text{ with }p(x)\text{ }SOS\text{,} \EA \]
and
\[\BA{c}
s\text{ is a moment sequence}\\
\Updownarrow \\
\text{ }\left\langle s,p_{\alpha }\right\rangle \geq 0,\quad
\text{for all }p(x)\in \reals_{m}[x]\text{ with }p(x)\geq 0\text{
on }\reals^{n}\text{,} \EA\] hence the cone of coefficients of$\
SOS$ polynomials and that of p.s.d. sequences are polar, and so
are the cones of moment sequences and positive polynomials.

From \cite{Puti93} then, we know that the problem of testing if a sequence $%
y $ is the moment sequence of some measure $\mu $ with support in
a compact semialgebraic set
\[
K=\left\{ x\in \reals^{n}:g_{k}(x)\geq 0,k=1,...,r\right\}
\]
and the problem of representing positive polynomials on $K$ as
$g_{k}(x)$ weighted sums of $SOS$ polynomials are dual of each
other and both representable as linear matrix inequalities.

\section{Semidefinite representation of the MinCard(x) and MinRank(X)
problems.}

\label{s:sdp-sap}As above $K$ is the semialgebraic set defined by
$K=\left\{ x\in \reals^{n}:g_{k}(x)\geq 0,k=1,...,r\right\} $ and
we assume that (\ref {assumpt-Putinar}) holds. Let $p(x)\in $
$\reals[x_{1},...,x_{n}]$ and as $K$ is compact, we also note
$t^{\ast }=\min_{K}p(x)$. Of course, $p(x)-t^{\ast
}\geq 0$ on $K$, hence there are $SOS$ polynomials $q_{k}(x)\in $ $%
\reals%
[x_{1},...,x_{n}]$ for $k=1,...,r$ such that
(\ref{prop-Putinar-sos}) holds for $p(x)-t^{\ast }$ on $K$. We
first show that the $\mathbf{Min}\Card(x)$ problem can be cast as
a semialgebraic program, hence a semidefinite program, using the
results from section \ref{s:sdr-sos}.

\begin{proposition}
Let $A\in \reals^{m\times n}$ and $b\in \reals^{m}$. There are polynomials $%
g_{k}(x)\in \reals[x_{1},...,x_{n}],$ for $k=0,...,r$ such that
the optimum values of:
\begin{equation}
\begin{array}{ll}
\text{minimize} & \Card(x) \\
\text{subject to} & Ax\geq b,
\end{array}
\label{pb-min-card-LMI}
\end{equation}
and
\[
\begin{array}{ll}
\text{minimize} & g_{0}(x) \\
\text{subject to} & g_{k}(x)\succeq 0,\quad \text{for }i=1,...,r,
\end{array}
\]
are equal.
\end{proposition}

\begin{proof}
First, as in \cite{Shor87}\ we notice that:
\[
\begin{array}{lll}
\Card(x)= & \min & \sum_{i=1}^{n}v_{i} \\
& \text{s.t.} & (v_{i}-1)x_{i}=0 \\
&  & v_{i}\geq 0,\quad \text{for }i=1,...,n,
\end{array}
\]
hence the $\mathbf{Min}\Card(x)$ problem in
(\ref{pb-min-card-LMI}) can be written:
\[
\begin{array}{lll}
\mathbf{Min}\Card(x)\equiv & \text{min.} & \sum_{i=1}^{n}v_{i} \\
& \text{s.t.} & (v_{i}-1)x_{i}=0 \\
&  & v_{i}\geq 0,\quad \text{for }i=1,...,n \\
&  & Ax\geq b.
\end{array}
\]
which is a semialgebraic problem.
\end{proof}

We now show a similar result on the $\mathbf{Min}\Rank(X)$, a
minimum cardinality problem on the eigenvalues of the matrix $X$.

\begin{proposition}
Let $A_{i}\in \symm^{n}$, for $i=1,...,p$ and $b\in \reals^{p}$.
There are polynomials $g_{k}(x)\in \reals[x_{1},...,x_{n}],$ for
$k=0,...,r$ such that the optimum values of:
\begin{equation}
\begin{array}{ll}
\text{minimize} & \Rank(X) \\
\text{subject to} & Tr\left( A_{i}X\right) =b_{i},\quad \text{for
}i=1,...,p
\\
& X\succeq 0
\end{array}
\end{equation}
and
\[
\begin{array}{ll}
\text{minimize} & g_{0}(x) \\
\text{subject to} & g_{k}(x)\succeq 0,\quad \text{for }i=1,...,M
\end{array}
\]
are equal.
\end{proposition}

\begin{proof}
We note $\left( \lambda _{1},...,\lambda _{n}\right) $ the
eigenvalues of the matrix $X\succeq 0$ and
\[
\sigma _{k}(X)=\sum_{\alpha \subset \{1,\ldots ,n\},\left| \alpha
\right| =k}\lambda _{\alpha }
\]
the symmetric functions. We note $\chi _{t}\left( X\right) $ the
characteristic polynomial of the matrix $X,$ with $\chi _{t}\left(
X\right) =\sum_{i=1}^{n}\left( -1\right) ^{i}\sigma _{i}(X)t^{i}$.
Because the matrix $X$ is semidefinite positive, we have $\sigma
_{k}(X)=0$ iff $\Rank\left( X\right) <k$, hence the
$\mathbf{Min}\Rank(X)$ problem can be expressed as a minimum
cardinality problem on the coefficients of the characteristic
polynomial. With
\[
\begin{array}{lll}
\Rank(X)= & \min & v_i \\
& \text{s.t.} & \sigma _{i}(X)(v_{i}-1)=0 \\
&             & v_i \geq 0,\quad \text{for }i=1,...,n,
\end{array}
\]
we enforce the remaining constraints to get the following
representation of $\mathbf{Min}\Rank(X)$:
\[
\begin{array}{ll}
\text{minimize} & \sum_{i=1}^{n}v_{i} \\
\text{s.t.} & (v_{i}-1)\sigma _{i}(X)=0 \\
& v_{i}\geq 0,\quad \text{for }i=1,...,n \\
& Tr\left( A_{i}X\right) =b_{i},\quad \text{for }i=1,...,p \\
& d_{I}(X)\geq 0,\quad \text{for }I\subset \{1,\ldots ,n\},
\end{array}
\]
where $d_{I}(X)$ is the principal minor with index set $I\subset
\{1,\ldots ,n\}$. This is a semialgebraic program in the
coefficients of the matrix $X$.
\end{proof}

These two results together with the results cited in section
\ref{s:sdr-sos} show that the two problems considered are
equivalent to very large scale semidefinite programs.

\section{Semidefinite relaxations}

\label{s:sdp-relax}In practice, the exact representations obtained
in the last section can be exponentially large and in general, we
cannot expect these problems to be tractable. Hence, the central
contribution of these representations is not to reduce the
complexity of these problems, but to provide a sequence of
successively sharper relaxations covering the entire complexity
spectrum, thus allowing the complexity/sharpness tradeoff to be
tuned. This is what we intend to describe in this section.

We begin by recalling the construction of moment matrices as
detailed in \cite{Curt00}, \cite{Lass01} and \cite{Lass02a}.
Again, we let \label{def-moment-vector}
$y_{x}=\left(1,x_{1},...,x_{n},x_{1}^{2},x_{1}x_{2},...,x_{1}x_{n},x_{2}x_{3},...,x_{n}^{2},...,x_{1}^{m},...,x_{n}^{m}\right)$
be the vector of all monomials in $\reals[x_{1},...,x_{n}]$, up to degree $m$%
, listed in increasing graded lexicographic order. We note $s(m)$
the size of the vector $y_{x}$. Let $y\in \reals^{s(2m)}$ be the
vector of moments (indexed according to $y_{x}$) of some
probability measure $\mu $ with
support $K =\left\{ x\in \reals^{n}:g(x)\geq 0\right\} $, we note $%
M_{m}(y)\in \symm^{s(m)}$, for the moment matrix defined by
\begin{equation}
M_{m}(y)_{i,j}=\int_{K}\left( y_{x}\right) _{i}\left( y_{x}\right)
_{j}\mu \left( dx\right) ,\quad \text{for }i,j=1,...,s(m)
\end{equation}
i.e. the (symmetric) matrix of moments with rows and columns indexed as in $%
y_{x}$. We note $\beta (i)$ the exponent of the monomial $\left(
y_{x}\right) _{i}$ and conversely, we note $i(\beta )$ the index
of the
monomial $x^{\beta }$ in $y_{x}$. For a given moment vector $y\in \reals%
^{s(m)}$ ordered as in (\ref{def-moment-vector}), the first row
and columns of the matrix $M_{m}(y)$ are then equal to $y$. The
rest of the matrix is then constructed following:
\[
M_{m}(y)_{i,j}=y_{\alpha +\beta }\text{ if
}M_{m}(y)_{1,i}=y_{\alpha }\text{ and }M_{m}(y)_{j,1}=y_{\beta }.
\]

Similarly, let $g(x)\in \reals[x_{1},...,x_{n}]$, we derive the
moment matrix for the measure $g(x)\mu \left( dx\right) $ on $K$
(called the localizing matrix), noted $M_{m}(gy)\in \symm^{s(m)}$,
from the matrix of moments $M_{m}(y)$ by:
\begin{equation}
M_{m}(gy)_{i,j}=\int_{K}\left( y_{x}\right) _{i}\left(
y_{x}\right) _{j}g(x)\mu \left( dx\right)
\end{equation}
for $i,j=1,...,s(m)$. The coefficients of the matrix $M_{m}(gy)$
are then given by:
\begin{equation}
M_{m}(gy)_{i,j}=\sum_{\alpha }g_{\alpha }M_{m}(y)_{i\left( \beta
(i)+\beta (j)+\alpha \right) }  \label{def-localizing-matrix}
\end{equation}
We can remark as in \cite{Lass01} that if the measure $\mu $ has
its support included in $K=\left\{ x\in \reals^{n}:g(x)\geq
0\right\} $, then for all coefficient vectors $v\in
\reals^{s(m)}$:
\[
\left\langle v,M_{m}(gy)v\right\rangle =\int_{K}v(x)^{2}g(x)\mu
\left( dx\right) \geq 0
\]
hence $M_{m}(gy)\succeq 0$.

In dimension one, for a given vector $y\in \reals^{s(2m)}$,
$M_{m}(y)\succeq 0$ (which is a LMI) is also a sufficient
condition in order for $y$ to the moment sequence of a probability
measure. In $\reals^{n}$, this equivalence does not hold in
general. The compact semialgebraic case is called the K-moment
problem and is dual to the compact $SOS$ problem in (\ref
{prop-Lasserre-sos}). Following \cite{Lass01}, we now exploit this
duality
to compute a sequence of semidefinite relaxations for the $\mathbf{Min}\Card(x)$%
and $\mathbf{Min}\Rank(X)$ problems.

\subsection{The MinCard(x) problem}

In section \ref{s:sdp-sap}, we saw that the optimum value of the
$\mathbf{Min}\Card(x)$ problem can be computed as the optimum
value of the semialgebraic program:
\begin{equation}
\begin{array}{ll}
\text{min.} & \sum_{i=1}^{n}v_{i} \\
\text{s.t.} & (v_{i}-1)x_{i}=0 \\
& v_{i}\geq 0,\quad \text{for }i=1,...,n \\
& a_{j}^{T}x\geq b_{j},\quad \text{for }j=1,...,m.
\end{array}
\end{equation}
As in \cite{Lass01}, to ensure compactness, we impose the
additional constraint $x_{1}^{2}+...+x_{n}^{2}\leq \alpha $ for
some constant $\alpha
>1 $. It is easy to check that the program above, together with this
additional bound on the feasible set, satisfies the constraints
qualification assumption (\ref{assumpt-Putinar}). For $N\geq 1$, a
lower bound $l_N$ on the optimal value of the above problem is
then computed as:
\begin{equation}
\begin{array}{lll}
l_N:= & \text{inf} & \sum_{i=1}^{n}y_{i} \\
& \text{s.t.} & M_{N}(y)\succeq 0 \\
&  & M_{N-1}\left( x_{i}(v_{i}-1)y\right) =0 \\
&  & M_{N-1}\left( (\alpha -x^{T}x-v^{T}v)y\right) \succeq 0 \\
&  & M_{N-1}\left( v_{i}y\right) \succeq 0,\quad \text{for }i=1,...,n \\
&  & M_{N-1}\left( \left( a_{j}^{T}x-b_{j}\right) y\right) \succeq
0,\quad \text{for }j=1,...,m,
\end{array}
\label{prog-mincard-sdp}
\end{equation}
in the variable $y\in \reals^{2n}$. Theorem 3.2 in \cite{Lass02a}
then states that there exists some $N^{\ast }$ such that
\[
l_N=\mathbf{Min}\Card(x),\quad \text{for all }N\geq N^{\ast },
\]
and the optimum is achieved whenever the rank of the matrices
$M_{N}(\left( ...\right) y)$ stabilizes.

\subsection{The MinRank(X) problem}

In section \ref{s:sdp-sap}, for $X\in \symm^{n}$, we saw that the
optimum of the $\mathbf{Min}\Rank(X)$ problem can be computed as
the optimum value of the semialgebraic program:
\[
\begin{array}{ll}
\text{min.} & \sum_{i=1}^{n}v_{i} \\
\text{s.t.} & (v_{i}-1)\sigma _{i}(X)=0 \\
& v_{i}\geq 0,\quad \text{for }i=1,...,n \\
& Tr\left( A_{j}X\right) =b_{j},\quad \text{for }j=1,...,p \\
& d_{I}(X)\geq 0,\quad \text{for }I\subset \lbrack 1,n],
\end{array}
\]
To further simplify this program, we can substitute to the $2^{n}$
constraints on the principal minors a more economical
semialgebraic constraint. The modified program then reads:
\[
\begin{array}{ll}
\text{min.} & \sum_{i=1}^{n}v_{i} \\
\text{s.t.} & (v_{i}-1)\sigma _{i}(X)=0 \\
& v_{i}\geq 0,\quad \text{for }i=1,...,n \\
& Tr\left( A_{i}X\right) =b_{i},\quad \text{for }i=1,...,p \\
& u^{T}Xu\geq 0,\quad \text{for }u\in \reals^{n},
\end{array}
\]
and again, to ensure compactness, we impose
$X^{T}X+v^{T}v+u^{T}u\leq \alpha
$ for some constant $\alpha >1$. If we set the variable $x=(u,X,v),$ for $%
N\geq \left\lceil \frac{n+1}{2}\right\rceil $, a lower bound $l_N$
on the optimal value of the above problem is computed as:
\begin{equation}
\begin{array}{lll}
l_N:= & \text{inf} & \sum_{i=1}^{n}u_{i} \\
& \text{s.t.} & M_{N}(y)\succeq 0 \\
&  & M_{N-\left\lceil \frac{i+1}{2}\right\rceil }\left(
(v_{i}-1)\sigma
_{i}(X)y\right)=0\\
&  & M_{N-1}\left( v_{i}y\right) \succeq 0\\
&  & M_{N-1}\left( \left( Tr\left( A_{i}X\right) -b_{i}\right)
y\right)=0\\
&  & M_{N-1}\left( \left( \alpha -X^{T}X+v^{T}v+u^{T}u\right)
y\right)
\succeq 0 \\
&  & M_{N-2}\left( \left( u^{T}Xu\right) y\right) \succeq 0
\end{array}
\label{prog-minrank-sdp}
\end{equation}
for $i=1,\ldots,n$ and $j=1,\ldots,p$, in the variable $y\in
\reals^{2n}$, where the matrices $M(q(x)y)$ are computed as in
(\ref{def-localizing-matrix}). Theorem 3.2 in \cite{Lass02a} then
states that there exists some $N^{\ast }$ such that
\[
l_N=\mathbf{Min}\Rank(X),\quad \text{for all }N\geq N^{\ast },
\]
and the optimum is reached whenever the rank of the matrices
$M_{N}(\left( ...\right) y)$ stabilizes. Alternatively, one could
use the fact that if we note $\chi _{t}(X)$ the characteristic
polynomial of $X$, then $X\succeq 0$ is equivalent to $\chi
_{-t^{2}}(X)$ being $SOS$ as a univariate polynomial in $t$.

\subsection{Convex envelope}

Suppose that instead of having only one $\mathbf{Min}\Card(x)$ or
$\mathbf{Min}\Rank(X)$ problem to solve, we need to solve a (long)
sequence of these problems with only some variation in the
constraints. Here, instead of computing an exact relaxation for
every instance of the problem, we are interested in finding an
efficient heuristic method for approximating the solution to all
the problems to be solved. The complexity of the first ''bound
design'' program will be high, but that of the subsequent programs
will then be much lower. The heuristics in \cite{Faze00} replaced the $\Card%
(x)$ (resp. $\Rank(X)$) functions by their convex envelope on the sets $%
0\leq x\leq 1$ (resp. $0\preceq X\preceq I$), i.e. the largest
convex
function $f(x)$ such that $f(x)\leq \Card(x)$ if $0\leq x\leq 1$ (resp. $%
f(X)\leq \Rank(X)$ if $0\preceq X\preceq I$). In this section, we
extend these bounds to semialgebraic sets with more complex
shapes.

Of course, a function and its convex envelope share the same
global minimum, so solving for this optimal lower bound is at
least as hard as finding the global minimum. Here however we look
for a convex lower bound for the problem in
(\ref{pb-min-card-LMI}) inside the set of polynomials of degree at
most $d$. This becomes a semialgebraic program:
\begin{equation}
\begin{array}{ll}
\text{maximize} & \int_{[0,1]^n} p(x)dx \\
\text{subject to} & t-p(x)\geq 0\text{ on }K \\
& \sum_{i=1}^{n}v_{i}-p(x)\geq 0\text{ on }K \\
& u^{T}\nabla ^{2}p(x)u\geq 0\text{ on }\left\| u\right\| ^{2}=1
\end{array}
\end{equation}
in the coefficients $p_{\alpha }$ of the polynomial $p(x)\in %
\reals_{d}[x_{1},...,x_{n}]$, where $K$ is the compact
semialgebraic set given by:
\begin{equation}
\begin{array}{ll}
K=(v,x)\in \reals^{2n}: & \left( v_{i}-1\right) x_{i}=0 \\
& Ax-b\geq 0 \\
& x,v\geq 0.
\end{array}
\end{equation}
Again, this can be cast as a LMI using the technique in
\cite{Lass02a}. We
notice that the $l_{1}$ heuristic is a particular case when the constraint $%
Ax-b\geq 0$ is dropped.

\section{Complexity}

\label{s:complexity}Of course, the two semidefinite programs
detailed in the last section are far from tractable if the
dimension $n$ and the relaxation order $N$ grow beyond textbook
example sizes. The $\mathbf{Min}\Card(x)$ problem is equivalent to
solving $2^{n}$ linear programs, so it is right to ask whether the
programs above provide any benefit over, for example,
branch-and-bound methods?

Even if these two methodologies have similar worst-case
complexities, the semidefinite relaxations in
(\ref{prog-mincard-sdp}) and (\ref{prog-minrank-sdp}) do sometimes
produce the global optimum for low order $N$ (see \cite{Lass02a})
and because the objective is integer valued here, they only need
to be solved up to an absolute precision of 1/2. We quickly detail
below some other possible simplifications.

But the results above have to be considered first as
representations, providing an insight of the relative complexity
of minimum cardinality problems versus that \ of tractable convex
optimization problems.

\subsection{Structure, sparsity and symmetry}

The first element that can be used to simplify the programs in
(\ref {prog-mincard-sdp}) and (\ref{prog-minrank-sdp}) is
\textit{structure}. In \cite{Lass02a} for example, the constraints
$x\in \{-1,1\}$, that translate into $x^{2}-x=0$ and $M_{m}(\left(
x^{2}-x\right) y)=0$ also imply that the variables $y_{\alpha },$
for $\alpha \in \mathbf{Z}_{+}^{n}$, can be replaced by
$y_{\Card(\alpha )}$ in the program. Secondly, if the constraints
in (\ref{prog-mincard-sdp}) and (\ref{prog-minrank-sdp}) include
some \textit{symmetry}, we could use the results in \cite{Gate00}
and \cite {Gate02} to preprocess and simplify the original
semialgebraic program. The simplest example of these symmetries is
of course when the constraints are invariant with respect to a
change of basis, in which case $\mathbf{Min}\Rank(X)$ reduces to a
$\mathbf{Min}\Card(x)$ problem and as in \cite{Gate02} p. 31, the
asymptotic complexity goes from being exponential in $n$ to being
exponential in $\sqrt{n}$.

In general, the complexity of algorithms in semialgebraic geometry
grows at least exponentially with the dimension. Efficiency is
then very often measured by the ability of a method to maintain
\textit{sparsity}. Some results on Newton polytopes and $SOS$
polynomials can then be used to efficiently handle sparsity in the
SOS representations. In particular, a
result of \cite{Rezn78} shows that if $p(x),h_{i}(x)\in \reals[x]$, for $%
i=1,...,r$, with $p(x)=\sum_{i=1}^{r}h_{i}^{2}(x)$, then
$C(h_{i})\subseteq \frac{1}{2}C(p)$. That result however does not
hold as is for the representation in (\ref{prop-Putinar-sos}).

Finally, a lower bound on the optimal value can be
obtained by simply dropping some of the constraints in (\ref{prog-mincard-sdp}%
) and (\ref{prog-minrank-sdp}).

\section{Conclusion}

One of the central contributions of semidefinite programs to the
optimization toolbox is their ability to efficiently solve a wide
class of convex eigenvalue problems. In this work, we have
illustrated how the method described in \cite{Lass01}\ for solving
semialgebraic programs, by lifting them to semidefinite programs,
can also be used to represent some semialgebraic eigenvalue
problems and convex envelope relaxations. This contribution is
centered around semidefinite representations and the insight they
can provide on the theoretical complexity of these problems.
Wether or not they also improve the practical complexity of
computing relaxations to these problems remains to be explored.

\bibliographystyle{amsalpha}
\providecommand{\bysame}{\leavevmode\hbox
to3em{\hrulefill}\thinspace}
\providecommand{\MR}{\relax\ifhmode\unskip\space\fi MR }
\providecommand{\MRhref}[2]{%
  \href{http://www.ams.org/mathscinet-getitem?mr=#1}{#2}
} \providecommand{\href}[2]{#2}

\end{document}